\documentclass[11pt]{article}
\usepackage{stmaryrd}
\usepackage{mathrsfs}
\usepackage{amsmath}
\usepackage{amsmath,amsthm,amssymb}
\usepackage{amsfonts}
\usepackage{latexsym}
\usepackage{CJK}
\date{}

\begin{document}
\title{Endomorphism algebras of semi-tilting modules$^\star$}
\author{{\small Shunhua Zhang}\\
         {\small School of Mathematics, Shandong University,
        Jinan, 250100,P.R.China}}

\pagenumbering{arabic}

\maketitle
\begin{center}
 \begin{minipage}{120mm}
   \small\rm
   {\bf  Abstract}\ \ Let $A$ be a finite dimensional algebra over an algebraically
   closed field $k$.  We investigate the structure properties of the endomorphism algebras of semi-tilting $A$-modules,  and prove that
   the endomorphism algebras arising from the mutations of semi-tilting $A$-modules can be realized as the endomorphism algebras of BB-tilting modules.

\end{minipage}
\end{center}

\vskip0.1in



\footnote{  Email addresses: \    shzhang@sdu.edu.cn.}

\vskip0.2in

\section {Introduction}

Tilting theory plays an important role in the representation theory of Artin algebras,
and the endomorphism algebras of tilting modules form a central
class of Artin algebras. The tilting modules of finite
projective dimension was defined by  Miyashita in \cite{M}. A further generalization of tilting modules
to modules of possibly infinite projective dimension called semi-tilting modules was made by Koga in \cite{K}.

\vskip0.2in

In \cite{K}, Koga proved that the class of basic semi-tilting modules is closed under mutations.
In this paper we investigate the structure properties of endomorphism algebras of semi-tilting modules
and prove that the endomorphism algebras of semi-tilting modules can be realised as the endomorphism algebras of BB-tilting modules.
We state our main results as follows.

\vskip 0.2in

{\bf Theorem 1.}\  {\it Let $A$ be a finite dimensional algebra over an algebraically
   closed field $k$, and $T=U\oplus X$ be a basic semi-tilting $A$-modules with $X$ indecomposable and $X\in {\rm gen}\ U$.
   Set  $B={\rm End}_{A}\ T$.  Then there exists an exact sequence
   $0\longrightarrow Y \stackrel{\mu}\longrightarrow E\stackrel{\varepsilon}\longrightarrow X\longrightarrow 0$
   in $A$-mod  with $\varepsilon$ a minimal right ${\rm add}\ U$-approximation,
   such that $T'=U\oplus Y$ is a  basic semi-tilting $A$-modules
   and $M={\rm Hom}_{A}(T', T)$ is a BB-tilting $B^{op}$-module with ${\rm End}_{B^{op}}\ M\simeq {\rm End}_{A}\ T'$.}

\vskip 0.2in

This paper is arranged as the following. In section 2, we fix the
notations and recall some necessary facts needed for our further
research. Section 3 is devoted to the proof of Theorem 1.

\vskip 0.2in

\section {Preliminaries}

\vskip 0.2in

Let $A$ be a finite dimensional algebra over an algebraically closed
field $k$. We denote by ${\rm gl.dim}\ A$ the global dimension of
$A$, and by $A$-mod the category of all finitely
generated left $A$-modules.  $D={\rm Hom}_k(-,\ k)$ is the
standard duality between $A$-mod and $A^{op}$-mod. $\tau$ is the
Auslander-Reiten translation of $A$ and  $\tau^{-1}$ is its inverse.

\vskip 0.2in

Given a $A$-module $M$, we denote by ${\rm add} \ M$ the full subcategory
having as objects the direct sums of indecomposable summands of $M$
and by ${\rm pd}_A M$ the projective dimension of $M$.
Let $M = M_1^{n_1}\oplus \cdots \oplus M_t^{n_t}$ where the $M_i$ are pairwise non-isomorphic indecomposable modules and
$n_i \geq 1$ is the multiplicity of $M_i$ in $M$. The module $M$ is called basic if $n_i = 1$  for all $i$.
We denote by $\delta (M)$ the number of non-isomorphic indecomposable summands of $M$.

\vskip 0.2in

  A module $T\in A$-mod is called a tilting module if
the following conditions are satisfied:\\
{(1)} ${\rm pd}_A T =n<\infty$;\\
{(2)} ${\rm Ext}_A^{i}(T,T)=0 $ for all $i >0$;\\
{(3)} There is  a long exact sequence
 $$0\longrightarrow
A\longrightarrow T_{0}\longrightarrow T_{1}
\longrightarrow\cdot\cdot\cdot\longrightarrow T_{n}\longrightarrow
0$$ with $T_{i}\in {\rm add}\ T $ for $0\leq i\leq n$.

\vskip 0.2in

Let $S$ be a non-injective simple $A$-module with the
following two properties:

(a) ${\rm pd}_A \tau^{-1}S  \leq 1$, and

(b) ${\rm Ext}^1_{A}(S, S) = 0$.

\vskip 0.2in

We denote the projective cover of $S$
by $P(S)$, and assume that $A = P(S) \oplus P$ such that there
is not any direct summand of $P$ isomorphic to $P(S)$. Let $T =
\tau^{-1} S \oplus P$ and $\Gamma={\rm End}_A \ T$. Then $T$ is a tilting
module with ${\rm pd}_A T \leq 1$, which is called BB-tilting module associated to $S$, and
$T_\Gamma$ is also a BB-tilting $\Gamma$-module, see Section 2.8 in \cite{A} for details.

\vskip 0.2in

Let $\mathcal{C}$ be a full subcategory of $A$-mod,
$C_{M}\in\mathcal{C}$ and $\varphi :C_M\longrightarrow M$ with
$M\in$ $A$-mod. Recall from \cite{AR} that the morphism $\varphi$ is a right
$\mathcal{C}$-approximation of $M$ if the induced  morphism ${\rm
Hom}_A(C,C_{M})\longrightarrow {\rm Hom}_A(C,M)$ is surjective for
any $C\in\mathcal{C}$. A minimal right $\mathcal{C}$-approximation
of $M$ is a right $\mathcal{C}$-approximation which is also a right
minimal morphism, i.e., its restriction to any nonzero summand is
nonzero. The subcategory $\mathcal{C}$ is called contravariantly
finite if any module $M\in$ $\Lambda$-mod admits a (minimal) right
$\mathcal{C}$-approximation. The notions of (minimal) left
$\mathcal{C}$-approximation and of covariantly finite subcategory
are dually defined. It is well known that add $M$ is both a
contravariantly finite subcategory and a covariantly finite
subcategory.

\vskip 0.2in

 The following Lemma is taken from \cite{K}.

\vskip 0.2 in

 {\bf Lemma 2.1.}  {\it Let $T=U\oplus X\in A$-mod with $X$ indecomposable,
 $X\not\in {\rm add} \ U$ and ${\rm Ext}^i_{A}(T,T) = 0$ for
 $i\geq 0$. Assume that there exists an exact sequence $0\rightarrow Y\rightarrow E
 \stackrel{\varepsilon}\longrightarrow X\rightarrow 0$ with $\varepsilon$ a right
 ${\rm add}\ U$-approximation. Set $B={\rm End}_{\Lambda} T$. Then ${\rm Hom}_{A}(U\oplus Y, T)$
 is a tilting $B$-module with ${\rm pd}_B {\rm Hom}_{A}(U\oplus Y, T)=1$.}

\vskip 0.2in

Recall from  \cite{K}, a module $T\in A$-mod is said to be a semi-tilting module if
the following conditions are satisfied:\\
{(i)} ${\rm Ext}_A^{i}(T,T)=0 $ for all $i >0$;\\
{(ii)} There is  a long exact sequence $0\rightarrow
A\rightarrow T_{0}\rightarrow T_{1}
\rightarrow\cdot\cdot\cdot\rightarrow T_{n}\rightarrow
0$ with $T_{i}\in {\rm add}\ T $ for $0\leq i\leq n$.

\vskip 0.2in

Throughout this paper, we follow the standard terminologies and
notations used in the representation theory of algebras, see \cite{ASS, ARS, HR, R}.

\vskip 0.2in

\section {Endomorphism algebras of semi-tilting modules}

In this section, we prove our main theorem.

\vskip 0.2in

{\bf Theorem 3.1.}\  {\it Let $A$ be a finite dimensional algebra over an algebraically
   closed field $k$, and $T=U\oplus X$ be a basic semi-tilting $A$-modules
   with $X$ indecomposable and $X\in {\rm gen}\ U$.
   Set  $B={\rm End}_{A}\ T$.  Then there exists an exact sequence
   $0\rightarrow Y \stackrel{\mu}\longrightarrow E\stackrel{\varepsilon}\longrightarrow X\rightarrow 0$
   in mod-$A$ with $\varepsilon$ a minimal right ${\rm add}\ U$-approximation,
   such that $T'=U\oplus Y$ is a  basic semi-tilting $A$-modules
   and $M={\rm Hom}_{A}(T', T)$ is a BB-tilting $B^{op}$-module with ${\rm End}_{B^{op}}\ M\simeq {\rm End}_{A}\ T'$.}

\vskip 0.1in

{\bf Proof.}\ \ According to Lemma 2.1, there exists an exact sequence
$$({\dag})\ \ \ \ \ \     0\rightarrow Y \stackrel{\mu}\longrightarrow E\stackrel{\varepsilon}\longrightarrow X\rightarrow 0$$
in mod-$A$ with $\varepsilon$ a minimal right ${\rm add}\ U$-approximation,
such that $T'=U\oplus Y$ is a  basic semi-tilting $A$-modules
and $M={\rm Hom}_{A}(T', T)$ is a tilting $B^{op}$-module with ${\rm pd}_B\ M=1$.

Note that $M= {\rm Hom}_{A}(U\oplus Y, T)= P_B\oplus L$ with $P_B$ is a projective $B^{op}$-module and $L={\rm Hom}_{A}(Y, T)$.

Applying ${\rm Hom}_{A}(-,T )$ to $(\dag)$ yields an exact sequence
$$
(\ddag)\ \ \ \ \ \ \  0\rightarrow {\rm
Hom}_{A}(X,T)\stackrel{\varepsilon^*} \longrightarrow {\rm
Hom}_{A}(E,T) \stackrel{\mu^*}\longrightarrow L \rightarrow 0 .
$$
Note that $(\ddag)$ is the minimal projective resolution of $L$.
Applying ${\rm Hom}_{B^{op}}(-, B)$ to $(\ddag)$ we get an exact
sequence of $B$-modules
$$
{\rm Hom}_{B^{op}}({\rm Hom}_{A}(E, T), B^{op}) \stackrel{\varepsilon^{**}}\longrightarrow
{\rm Hom}_{B^{op}}({\rm Hom}_{A}(X,T),B^{op}) \longrightarrow
{\rm Tr}_{B^{op}}\ L \rightarrow 0,
$$
which is isomorphic to the following exact sequence
$$
{\rm Hom}_{A}(T,E)\stackrel{\varepsilon_*} \longrightarrow {\rm
Hom}_{A}(T,X) \longrightarrow {\rm Tr}_{B^{op}}\ L\rightarrow
0,
$$
where $\varepsilon_*={\rm Hom}_{A}(T,\varepsilon)$.

We claim that ${\rm Im}\ \varepsilon_*$ is the radical of the indecomposable
projective $B^{op}$ module ${\rm Hom}_{A}(T,X)$.

Indeed,  $\varepsilon_*={\rm Hom}_{A}(U\oplus X, \varepsilon)={\rm
Hom}_{A}(U, \varepsilon)\oplus {\rm Hom}_{A}(X, \varepsilon)$.

Since $\varepsilon$ is a minimal left ${\rm add}  U$-approximation of $X$, by
applying ${\rm Hom}_{A}(U, -)$ to $(\dag)$ we get an exact
sequence
$$
0\rightarrow {\rm Hom}_{A}(U,Y) \longrightarrow {\rm
Hom}_{A}(U, E) \stackrel{{\rm Hom}_{A}(U,
\varepsilon)}\longrightarrow {\rm Hom}_{A}(U, X)\rightarrow 0.
$$
Hence ${\rm Hom}_{A}(U, \varepsilon)$ is surjective.

By applying ${\rm Hom}_{A}(X, -)$ to $(\dag)$ we have an exact
sequence
$$
{\rm Hom}_{A}(X, E) \stackrel{{\rm Hom}_{A}(X, \varepsilon)}\longrightarrow {\rm Hom}_{A}(X,X) \longrightarrow {\rm
Ext}_{A}^1(X, Y)\rightarrow 0,
$$
it forces that ${\rm Im~Hom}_{A}(X, \varepsilon)= {\rm rad}\ {\rm
Hom}_{A}(X,X)$ since ${\rm dim}_k\ {\rm Ext}_{A}^1(X, Y)=1$. It
follows that ${\rm Im}\ \varepsilon_*={\rm rad}\ {\rm Hom}_{A}(T, X)$, and our
claim is true.

It follows from our claim that ${\rm Tr}_{B^{op}}\ L$ is a simple
$B^{op}$-module and $\tau_{B^{op}} L$ is the simple socle $S$ of the
indecomposable injective $B^{op}$-module $D{\rm
Hom}_{A}(T, X)$, and we know that $L\simeq \tau_{B^{op}}^{-1}
S$. Hence $M$ is a BB-tilting $B^{op}$-module. By Lemma 3.4 in
\cite{HX}, we have that ${\rm End}_{B^{op}}\ M =
{\rm Hom}_{B^{op}}  ({\rm Hom}_{A}(U\oplus Y, T), {\rm Hom}_{A}(U\oplus Y, T))\simeq {\rm
End}_{A}\ (U\oplus Y)$. $\hfill\Box$

\vskip 0.2in

We state the dual result as follows.

\vskip 0.2in

{\bf Theorem 3.2.}\  {\it Let $A$ be a finite dimensional algebra over an algebraically
   closed field $k$, and $T=U\oplus Y$ be a basic semi-tilting $A$-modules with $Y$ indecomposable and $Y\in {\rm cogen}\ U$.
   Set  $B={\rm End}_{A}\ T$.  Then there exists an exact sequence
   $0\rightarrow Y \stackrel{\mu}\longrightarrow E\stackrel{\varepsilon}\longrightarrow X\rightarrow 0$
   in mod-$A$ with $\mu$ a minimal left ${\rm add}\ U$-approximation, such that $T'=U\oplus X$ is a  basic semi-tilting $A$-modules
   and $M={\rm Hom}_{A}(T', T)$ is a BB-tilting $B^{op}$-module with ${\rm End}_{B^{op}}\ M\simeq {\rm End}_{A}\ T'$.}

\vskip 0.2in

{\bf Remark.}\ Let $T=U\oplus X$ be a basic semi-tilting $A$-modules
   with $X$ indecomposable and $X\in {\rm gen}\ U$. We denote by $\mu_X (T)$ the module $T'=U\oplus Y$ in the case of Theorem 3.1.

\vskip 0.2in

Recall from \cite{K}, the semi-tilting quiver $\mathcal{K}$ is defined as follows: The vertices of $\mathcal{K}$ are
isomorphism classes of basic semi-tilting modules and there is an arrow $V \rightarrow W$ if
$W$ and $V$ are represented by basic semi-tilting modules $T$ and $\mu_X(T)$ with $X$ a
non-projective indecomposable direct summand of $T$, respectively.

\vskip 0.2in

{\bf Corollary 3.3.} \  {\it Let $A$ be a finite dimensional algebra over an algebraically
   closed field $k$, and $\mathcal{C}$ be a connected component of the semi-tilting quiver $\mathcal{K}$ of $A$.
   Then all the endomorphism algebras of semi-tilting modules in $\mathcal{C}$ can be realized
   as the iterated endomorphism algebras of BB-tilting modules.}

\vskip 0.1in

{\bf Proof}\ \  Let $T_1$ and $T_2$ be a pair of basic semi-tilting modules in $\mathcal{C}$.
Then there is a path
$T_1= W_0\rightarrow W_1\rightarrow\cdots\rightarrow W_s\rightarrow T_2=W_{s+1}$ in $\mathcal{C}$.
Set $\Lambda_i= {\rm End}_{A} W_i$ for $0\leq i\leq s+1$.
According to Theorem 3.2, for each arrow
$W_i\rightarrow W_{i+1}$ there is a BB-tilting $\Lambda_i^{op}$-module $M_i$
such that ${\rm End}_{\Lambda_i^{op}}\ M_i= \Lambda_{i+1} $.
In particular,  ${\rm End}_{A}\ T_2\simeq  {\rm End}_{\Lambda_{s}^{op}}\ M_{s}$.
This completes the proof.     $\hfill\Box$

\vskip 0.2in

Note that a semi-tilting module is a tilting module if and only if its projective dimension is finite.
According to \cite{HU2}, the tilting quiver $\mathcal{T}$ of $A$ is defined as follows. The vertices of $\mathcal{T}$ are
isomorphism classes of basic tilting modules,  and there is an arrow $T' \rightarrow T$ if
$T'=\overline{T}\oplus X$ and $T=\overline{T}\oplus Y$  with $X,
Y$ indecomposable and there is a short exact sequence
$0\rightarrow X\stackrel{f}\longrightarrow E\stackrel{g}\longrightarrow Y\rightarrow 0$
with $f$ (resp.$g$) being a minimal left (resp.right)
${\rm add}\ \overline{T}$-approximation.

\vskip 0.2in

Let $\mathcal{C}$ be a connected component of be the semi-tilting quiver $\mathcal{K}$ of $A$.
According to Theorem 3.11 in \cite{K} we know that $\mathcal{C}$ is a connected component of  the  tilting quiver $\mathcal{T}$ of $A$
if and only if $\mathcal{C}$ containing a tilting $A$-module. The following result is a consequence of Corollary 3.3.

\vskip 0.2in

{\bf Corollary 3.4.} \  {\it Let $A$ be a finite dimensional algebra over an algebraically
   closed field $k$, and $\mathcal{C}$ be a connected component of the tilting quiver $\mathcal{T}$ of $A$.
   Then all the endomorphism algebras of tilting modules in $\mathcal{C}$ can be realized
   as the iterated endomorphism algebras of BB-tilting modules.}

\end{document}